\documentclass[11pt,reqno]{amsart}
\usepackage{amssymb,amsmath,amsthm,mathtools,amsfonts,times,hyperref,mathrsfs,multirow,tabularx}
\usepackage{pgfplots}
\pgfplotsset{compat=1.15}
\usepackage{mathrsfs,mathdots}
\usepackage{esvect}
\usetikzlibrary{arrows}
\usepackage{enumerate}
\usepackage{graphicx}
\usepackage{array}
\usepackage{ytableau}
\usepackage{pstricks,pst-node,graphicx}
\usepackage{tikz,varwidth}
\usepackage{setspace}
\usepackage{float}
\usetikzlibrary{calc}
\usepackage[super]{nth}
\usepackage{tikz}
\usepackage{extarrows}
\usepackage{pgfplots}
\pgfplotsset{compat=1.15}
\usetikzlibrary{arrows}

\newtheorem{theorem}{Theorem}[section]

\theoremstyle{definition}

\theoremstyle{remark}

\numberwithin{equation}{section}

\allowdisplaybreaks

\newcommand\numberthis{\addtocounter{equation}{1}\tag{\theequation}}
\newcommand\restr[2]{{
  \left.\kern-\nulldelimiterspace 
  #1 
  \littletaller 
  \right|_{#2} 
  }}
\newcommand{\littletaller}{\mathchoice{\vphantom{\big|}}{}{}{}}  

\begin{document}

\title[$\vv{M}$ Versions of Andrews-Gordon Identities Revisited]{$\vv{M}$ Versions of Andrews-Gordon Identities Revisited}

\author{Alexander Berkovich}
\address{Department of Mathematics, University of Florida, Gainesville
FL 32611, USA}
\email{alexb@ufl.edu}
\author{Aritram Dhar}
\address{Department of Mathematics, University of Florida, Gainesville
FL 32611, USA}
\email{aritramdhar@gmail.com}

\date{\today}

\subjclass[2020]{05A17, 05A19, 11P81, 11P84}             

\keywords{Andrews-Gordon identities, Berkovich-Paule variants}

\begin{abstract}
In this paper, we revisit the work of Berkovich and Paule on variants of the Andrews-Gordon identities. We find a generalization of their principal polynomial identity and, as a consequence, obtain new $\Vec{M}$ versions of the Andrews-Gordon identities. More precisely, for non-negative integers $M_1\ge M_2\ge M_3\ge\ldots\ge M_\nu$, we show that a broad class of multi-sums of the form
$$
\sum\limits_{\mathbf n}
\frac{
q^{
N_1^2+\cdots+N_\nu^2
-
M_1N_1
-
M_2N_2
-
\cdots
-
M_\nu N_\nu
}
}{
(q)_{n_1}(q)_{n_2}\cdots(q)_{n_\nu}
}
$$
can be expressed as a sum of products. Above, we use standard notations for $q$-Pochhammer symbols and $N_i = \sum\limits_{k=i}^{\nu}n_k$ for $1\le i\le \nu$. 
\end{abstract}

\maketitle

\section{Introduction}\label{s1}
In 1961, Gordon \cite{Gordon} established a natural generalization of the classical Rogers-Ramanujan partition theorem. We recall Gordon's theorem in the following form.

\begin{theorem}[Gordon \cite{Gordon}]\label{thm:Gordon}
Let $\nu\geq 1$ and $1\leq s\leq\nu+1$. Then the number of partitions of $N$ whose frequency representation $N = \sum_{j\geq 1}jf_j$ satisfies $f_1\leq s-1$ and $f_j + f_{j+1}\leq\nu$, $f_j\geq 0$ $(j\geq 1)$ is equal to the number of partitions of $N$ into parts that are not congruent to $0$ or $\pm s\pmod{2\nu+3}$.
\end{theorem}

Andrews \cite{Andrews} subsequently obtained an analytic
counterpart of Gordon's Theorem \ref{thm:Gordon}. His result can be expressed as the following identity.

\begin{theorem}[Andrews \cite{Andrews}]\label{thm:Andrews}
Let $\nu\geq 1$ and $1\leq s\leq\nu+1$. For $|q|<1$, we have
\begin{align*}
\sum\limits_{n_1,n_2,\ldots,n_\nu\geq 0}\frac{q^{N_1^2 + \cdots + N_\nu^2 + N_s + \cdots + N_\nu}}{(q)_{n_1}(q)_{n_2}\cdots(q)_{n_\nu}} &= \frac{1}{(q)_{\infty}}\sum\limits_{j = -\infty}^{\infty}(-1)^jq^{\frac{j((2\nu+3)(j+1) - 2s)}{2}}\\
&= \prod\limits_{\substack{j\geq 1\\j\not\equiv 0,\pm s\pmod{2\nu+3}}}\frac{1}{1-q^j},\numberthis\label{eq:andrews-gordon}
\end{align*}
where $$N_i = \begin{cases}n_i + n_{i+1} + \cdots + n_\nu, & \text{if}\,\,\, 1\leq i\leq\nu, \\ 0, & \text{if}\,\,\, i=\nu+1.\end{cases}$$
\end{theorem}
We employ the standard $q$-Pochhammer notation above and throughout the rest of the paper. For any complex variable $a$, we define $$(a;q)_{\infty} = (a)_{\infty} := \prod\limits_{j\geq0}(1-aq^j),$$ and, for any non-negative integer $m$, $$(a;q)_m = (a)_m := \prod\limits_{j\geq0}\frac{1-aq^j}{1-aq^{j+m}}.$$ The product on the right-hand side of \eqref{eq:andrews-gordon} is the generating function for partitions into parts that are not congruent to $0$ or $\pm s$ modulo $2\nu+3$. Thus, Andrews' identity provides an analytic interpretation of Gordon's partition theorem: the multi-sum on the left-hand side encodes the frequency restrictions appearing in Gordon's theorem, while the infinite product records the corresponding congruence restrictions on the parts.

We also define the $q$-binomial (Gaussian) coefficient as
\begin{align*}
{m\brack n}_q := \Bigg\{\begin{array}{lr}
\dfrac{(q)_m}{(q)_n(q)_{m-n}}\quad\text{for } m\ge n\ge 0,\\
0\qquad\qquad\quad\text{otherwise}.\end{array}
\end{align*}

Define $$\alpha_{i,s} := \max(0,i-s).$$

Next, for $L\equiv s+b\pmod{2}$, we introduce the polynomials $B_{s,b}^{\nu}(L,q)$ by
\begin{align*}
B_{s,b}^{\nu}(L,q) = B_{s,b}^{\nu}(L) &:= \sum\limits_{j=-\infty}^{\infty}\Bigg\{q^{j((2j+1)(2\nu+3)-2s)}{L\brack (2\nu+3)j+\frac{L+b-s}{2}}_q\\
&\quad - q^{(2j+1)((2\nu+3)j+s)}{L\brack (2\nu+3)j+\frac{L+b+s}{2}}_q\Bigg\}.\numberthis\label{eq:Bdef}    
\end{align*}

Using recurrences for $q$-binomial coefficients, it is easy to check that
\begin{align*}
B_{s,b}^{\nu}(L) = B_{s,b-1}^{\nu}(L-1) + B_{s,b+1}^{\nu}(L-1) + (q^{L-1}-1)B_{s,b}^{\nu}(L-2).\numberthis\label{eq:recurrenceB}    
\end{align*}
Also, the following symmetry relations hold.
\begin{align*}
B_{s,\nu+2}^{\nu}(L) &= B_{2\nu+3-s,\nu+1}^{\nu}(L),\numberthis\label{eq:rel1B}\\
B_{s,1}^{\nu}(L) &= B_{s,2}^{\nu}(L-1).\numberthis\label{eq:rel2B}
\end{align*}

Foda and Quano obtained the following finite analogue of Theorem \ref{thm:Andrews} in  \cite{FodaQuano}.
\begin{theorem}[Foda and Quano \cite{FodaQuano}]\label{thm:FodaQuano}
We have
\begin{align*}
\sum\limits_{\mathbf n}&q^{N_1^2+\cdots+N_\nu^2+N_s+\cdots+N_\nu}\prod\limits_{i=1}^{\nu}{n_i+L-2\sum\limits_{j=1}^{i}N_j-\alpha_{i,s-1}\brack n_i}_q\\
&= \Bigg\{\begin{array}{lr}
B_{s,\nu+1}^{\nu}(L),\qquad\quad\text{if } L\not\equiv s+\nu\pmod{2},\\
B_{2\nu+3-s,\nu+1}^{\nu}(L),\,\,\,\text{otherwise}.\end{array}
\end{align*}
\end{theorem}

It is straightforward to check that in the limit $L\rightarrow\infty$, Theorem \ref{thm:FodaQuano} reduces to Theorem \ref{thm:Andrews}. We next recall a $q$-multinomial coefficient introduced in this context as
\begin{align*}
{L\brack a}^{(\nu)} := \sum\limits_{j_1+j_2+\cdots+j_\nu = a+\frac{\nu L}{2}}q^{\sum\limits_{r=2}^{\nu}j_{r-1}(L-j_r)}{L\brack j_\nu}_q{j_\nu\brack j_{\nu-1}}_q\cdots{j_2\brack j_1}_q.    
\end{align*}

In \cite{Warnaar}, Warnaar subsequently proved the following finite identity by a recursive argument.
\begin{theorem}[Warnaar \cite{Warnaar}]\label{thm:Warnaar}
For $1\leq s\leq \nu+1$, we have
\begin{align*}
\sum\limits_{\mathbf n}&q^{N_1^2+\cdots+N_\nu^2+N_s+\cdots+N_\nu}\prod\limits_{i=1}^{\nu}{n_i+iL-2\sum\limits_{j=1}^{i}N_j-\alpha_{i,s-1}\brack n_i}_q\\
&= \Bigg\{\begin{array}{lr}
W_{s}^{\nu}(L,q),\qquad\quad\text{if } \nu(L+1)+s\equiv 1\pmod{2},\\
W_{2\nu+3-s}^{\nu}(L,q),\,\,\,\,\text{otherwise}.\end{array}\numberthis\label{eq:Warnaar}
\end{align*}
where
\begin{align*}
W_{s}^{\nu}(L,q) := \sum\limits_{j=-\infty}^{\infty}&\Bigg\{q^{j((2j+1)(2\nu+3)-2s)}{L\brack (2\nu+3)j+\frac{\nu+1-s}{2}}^{(\nu)}\\
&- q^{(2j+1)((2\nu+3)j+s)}{L\brack (2\nu+3)j+\frac{\nu+1+s}{2}}^{(\nu)}\Bigg\}.\numberthis\label{eq:Wdef}
\end{align*}
\end{theorem}

The preceding results provide various finite analogs of the Andrews-Gordon identities and, in particular, motivate the study of other types of generalizations.

The remainder of the paper is organized as follows. In Section \ref{s2}, we recall the Schilling-Warnaar finitization of the Andrews-Gordon identities and introduce the associated $q$-supernomial coefficients. These coefficients provide the principal algebraic framework used in the subsequent sections. In Section \ref{s3}, we revisit an identity of Berkovich and Paule and establish its generalization to the setting of $\Vec{M}$ versions. We also record a special case and two examples of the resulting identity. Finally, in Section \ref{s4}, we state new $\Vec{M}$ versions of the Andrews-Gordon identities and illustrate them with three explicit examples.

\section{Schilling-Warnaar Finitization of Andrews-Gordon Identities and $q$-Supernomial Coefficients}\label{s2}
Let $\Vec{L} = (l_1,l_2,\ldots,l_\nu)$ and define $L_i := l_i+l_{i+1}+\cdots+l_\nu$ for $1\leq i\leq\nu$. We begin by recalling the finite analog of the Andrews-Gordon identities due to Schilling and Warnaar \cite{SchillingWarnaar}.
\begin{theorem}[Schilling-Warnaar \cite{SchillingWarnaar}]\label{thm:SchillingWarnaar}
For $1\leq s\leq\nu+1$, we have
\begin{align*}
&\sum_{\mathbf n}
q^{N_1^2+\cdots+N_\nu^2+N_s+\cdots+N_\nu}
\prod_{i=1}^{\nu}
\left[
\begin{matrix}
n_i+\displaystyle\sum_{j=1}^{i}L_j
-2\displaystyle\sum_{j=1}^{i}N_j-\alpha_{i,s-1}\\
n_i
\end{matrix}
\right]_q
\nonumber\\
&=
\begin{cases}
\begin{aligned}
&\displaystyle
\sum_{j=-\infty}^{\infty}
\Biggl(
q^{j((2j+1)(2\nu+3)-2s)}
\left[
\begin{matrix}
\vec L\\
(2\nu+3)j+\dfrac{\nu+1-s}{2}
\end{matrix}
\right]_q
\\[-2pt]
&\displaystyle\qquad
-q^{(2j+1)((2\nu+3)j+s)}
\left[
\begin{matrix}
\vec L\\
(2\nu+3)j+\dfrac{\nu+1+s}{2}
\end{matrix}
\right]_q
\Biggr)
\end{aligned},
&
\displaystyle
\text{if }\sum_{j=1}^{\nu}L_j+s+\nu\equiv1\pmod{2},
\\[16pt]
\begin{aligned}
&\displaystyle
\sum_{j=-\infty}^{\infty}
\Biggl(
q^{j((2j+1)(2\nu+3)-2s)}
\left[
\begin{matrix}
\vec L\\
(2\nu+3)j+\dfrac{\nu+2-s}{2}
\end{matrix}
\right]_q
\\[-2pt]
&\displaystyle\qquad
-q^{(2j+1)((2\nu+3)j+s)}
\left[
\begin{matrix}
\vec L\\
(2\nu+3)j+\dfrac{\nu+2+s}{2}
\end{matrix}
\right]_q
\Biggr)
\end{aligned},
&
\displaystyle
\text{if }\sum_{j=1}^{\nu}L_j+s+\nu\equiv0\pmod{2},
\end{cases}\numberthis\label{eq:SchillingWarnaar}
\end{align*}
where
\[
\left[
\begin{matrix}
\vec L\\
a
\end{matrix}
\right]_q
=
\sum_{\substack{j_1,\ldots,j_\nu\\
j_1+\cdots+j_\nu
=
a+\frac12(l_1+2l_2+3l_3+\cdots+\nu l_\nu)}}
q^{\displaystyle\sum_{r=2}^{\nu}j_{r-1}(L_r-j_r)}
\left[
\begin{matrix}
l_\nu\\
j_\nu
\end{matrix}
\right]_q
\left[
\begin{matrix}
l_{\nu-1}+j_\nu\\
j_{\nu-1}
\end{matrix}
\right]_q
\cdots
\left[
\begin{matrix}
l_1+j_2\\
j_1
\end{matrix}
\right]_q
\]
is known as the $q$-supernomial coefficient.
\end{theorem}

It is immediate that \eqref{eq:SchillingWarnaar} reduces to \eqref{eq:andrews-gordon} as $l_\nu\rightarrow\infty$.

\section{The Berkovich-Paule Identity and its Generalization}\label{s3}
We now turn to the identity of Berkovich and Paule that motivates our main result. Their identity expresses a finite Andrews-Gordon type multi-sum as a sum of products of polynomials.
\begin{theorem}[Berkovich-Paule \cite{BerkovichPaule1,BerkovichPaule2}]\label{thm:BerkovichPaule1}
We have
\begin{align*}
\sum\limits_{\mathbf n}q^{N_1^2+\cdots+N_\nu^2-MN_1}\prod\limits_{i=1}^{\nu}{n_i+L-2\sum\limits_{j=1}^{i}N_j\brack n_i}_q = \sum\limits_{s=1}^{\nu+1}\widetilde{F}_s^{\nu}(M,1/q)\widetilde{F}_s^{\nu}(L-M,q),\numberthis\label{eq:BerkovichPaule1}   
\end{align*}
where
\begin{align*}
\widetilde{F}_s^{\nu}(L,q) := \Bigg\{\begin{array}{lr}
B_{s,\nu+1}^{\nu}(L,q),\qquad\quad\text{if } L+s+\nu\equiv 1\pmod{2},\\
B_{2\nu+3-s,\nu+1}^{\nu}(L,q),\,\,\,\,\text{otherwise}.\end{array}\numberthis\label{eq:Ftilde}    
\end{align*}
\end{theorem}

Now, we state a generalization of \eqref{eq:BerkovichPaule1}. Our goal is to extend this identity by allowing independent parameters. To this end, let $M_1\geq M_2\geq\cdots\geq M_\nu$ be non-negative integers. For $1\le s\le \nu+1$ and $\boldsymbol{l} := (l_1,l_2,\ldots,l_\nu)$, define
\[
G_s(\boldsymbol{l},q) = G_s(l_1,l_2,\ldots,l_\nu,q)
:=
\sum\limits_{\textbf{n}}
q^{N_1^2+\cdots+N_v^2+N_s+\cdots+N_\nu}
\prod\limits_{i=1}^{\nu}
\begin{bmatrix}
n_i+\displaystyle\sum_{j=1}^{i}L_j
-2\displaystyle\sum_{j=1}^{i}N_j-\alpha_{i,s-1}\\[2mm]
n_i
\end{bmatrix}_q,
\]
where $N_{\nu+1}=0$ and $L_j := l_j+l_{j+1}+\cdots+l_\nu$ for $1\leq j\leq\nu$. Then, we have
\begin{align*}
\sum\limits_{\mathbf n}
q^{N_1^2+\cdots+N_\nu^2-M_1N_1-\cdots-M_\nu N_\nu}
\prod\limits_{i=1}^{\nu}
\begin{bmatrix}
n_i+\displaystyle\sum_{j=1}^{i}L_j
-2\displaystyle\sum_{j=1}^{i}N_j\\[2mm]
n_i
\end{bmatrix}_q
=
\sum\limits_{s'=1}^{\nu+1}
G_{s'}(\boldsymbol{l}-\widetilde{\mathbf M},q)\,
G_{s'}\!\left(\widetilde{\mathbf M},1/q\right),\numberthis\label{eq:mainresult1}    
\end{align*}
where
\[
\widetilde{\mathbf M}
=
(\widetilde M_1,\widetilde M_2,\ldots,\widetilde M_\nu),
\]
with
\[
\widetilde M_i=
\begin{cases}
M_i-M_{i+1}, & 1\le i\le \nu-1,\\[1mm]
M_v, & i=\nu,
\end{cases}
\]
provided the following inequalities hold
\[
l_1\ge M_1-M_2\ge 0,
\]
\[
l_2\ge M_2-M_3\ge 0,
\]
\[
\vdots
\]
\[
l_{\nu-1}\ge M_{\nu-1}-M_\nu\ge 0,
\]
\[
l_\nu\ge M_\nu\ge 0.
\]

We now state another identity of Berkovich and Paule (see \cite[eq. (4.6)]{BerkovichPaule2}) which turns out to be the special case $\boldsymbol{l} = (0,\ldots,0,L)$ and $M_i = M$ for $1\leq i\leq \nu$ of \eqref{eq:mainresult1}. It can be stated as follows.
\begin{theorem}[Berkovich-Paule \cite{BerkovichPaule2}]\label{thm:BerkovichPaule2}
We have
\begin{align*}
\sum\limits_{\mathbf n}&q^{N_1^2+\cdots+N_\nu^2-M(N_1+\cdots+N_\nu)}\prod\limits_{i=1}^{\nu}{n_i+iL-2\sum\limits_{j=1}^{i}N_j\brack n_i}_q\\
&\quad = \sum\limits_{s'=1}^{\nu+1}G_{s'}(0,\ldots,0,L-M,q)G_{s'}(0,\ldots,0,M,1/q).\numberthis\label{eq:BerkovichPaule2}
\end{align*}
\end{theorem}

The proof of \eqref{eq:mainresult1} will be given elsewhere. We now state two examples to illustrate \eqref{eq:mainresult1}.

\subsection{Example 1}
Consider $\nu=2$ with
\[
l_1=4,\qquad l_2=3,\qquad M_1=3,\qquad M_2=2.
\]
Then
\[
M_1-M_2=1,\qquad M_2=2,
\]
and hence
\[
\widetilde{\mathbf M}=(1,2),
\qquad
\boldsymbol{l}-\widetilde{\mathbf M}=(3,1).
\]
Moreover, the required inequalities are satisfied since
\[
l_1=4\geq M_1-M_2=1,
\qquad
l_2=3\geq M_2=2.
\]

Since
\[
N_1=n_1+n_2,\qquad N_2=n_2,
\]
we have
\[
N_1^2+N_2^2-3N_1-2N_2
=n_1^2+2n_1n_2+2n_2^2-3n_1-5n_2.
\]
Thus, \eqref{eq:mainresult1} gives
\[
\begin{aligned}
&\sum_{n_1,n_2\geq0}
q^{n_1^2+2n_1n_2+2n_2^2-3n_1-5n_2}
\left[\begin{matrix}
n_1+7-2N_1\\
n_1
\end{matrix}\right]_q
\left[\begin{matrix}
n_2+10-2N_1-2N_2\\
n_2
\end{matrix}\right]_q
\\
&\hspace{2cm}
=
\sum_{s'=1}^{3}
G_{s'}((3,1),q)\,
G_{s'}((1,2),1/q),
\end{aligned}
\]
where
\[
G_s((3,1),q)
=
\sum_{n_1,n_2\geq0}
q^{N_1^2+N_2^2+N_s+\cdots+N_2}
\left[\begin{matrix}
n_1+4-2N_1-\alpha_{1,s-1}\\
n_1
\end{matrix}\right]_q
\left[\begin{matrix}
n_2+5-2N_1-2N_2-\alpha_{2,s-1}\\
n_2
\end{matrix}\right]_q
\]
and
\[
G_s((1,2),q)
=
\sum_{n_1,n_2\geq0}
q^{N_1^2+N_2^2+N_s+\cdots+N_2}
\left[\begin{matrix}
n_1+3-2N_1-\alpha_{1,s-1}\\
n_1
\end{matrix}\right]_q
\left[\begin{matrix}
n_2+5-2N_1-2N_2-\alpha_{2,s-1}\\
n_2
\end{matrix}\right]_q.
\]

Consequently, the identity becomes
\begin{align*}
\sum\limits_{n_1,n_2\geq0}q^{n_1^2+2n_1n_2+2n_2^2-3n_1-5n_2}&\left[\begin{matrix}n_1+7-2N_1\\n_1\end{matrix}\right]_q\left[\begin{matrix}n_2+10-2N_1-2N_2\\n_2\end{matrix}\right]_q\\
&= \sum_{s'=1}^{3}G_{s'}((3,1),q)G_{s'}((1,2),1/q).\numberthis\label{eq:example13.3'}
\end{align*}

Substituting the relevant parameters into the definition of $G_s$ gives the six polynomials appearing on the right-hand side of \eqref{eq:example13.3'}. Evaluating these polynomials explicitly and simplifying yields
\begin{align*}
\sum\limits_{n_1,n_2\geq0}q^{n_1^2+2n_1n_2+2n_2^2-3n_1-5n_2}&\left[\begin{matrix}n_1+7-2N_1\\n_1\end{matrix}\right]_q\left[\begin{matrix}n_2+10-2N_1-2N_2\\n_2\end{matrix}\right]_q\\
&= 2q^{-3}+6q^{-2}+8q^{-1}+15+14q+11q^2+7q^3+2q^4.\numberthis\label{eq:example13.3}
\end{align*}

Maple verifies \eqref{eq:example13.3}.

\subsection{Example 2}

Consider $\nu=3$ with
\[
l_1=3,\qquad l_2=3,\qquad l_3=4,\qquad
M_1=3,\qquad M_2=2,\qquad M_3=1.
\]
Then
\[
M_1-M_2=M_2-M_3=M_3=1,
\]
and hence
\[
\widetilde{\mathbf M}=(1,1,1),
\qquad
\boldsymbol{l}-\widetilde{\mathbf M}=(2,2,3).
\]
Moreover, the required inequalities are satisfied since
\[
l_1=3\geq M_1-M_2=1,\qquad
l_2=3\geq M_2-M_3=1,\qquad
l_3=4\geq M_3=1.
\]

Since
\[
N_1=n_1+n_2+n_3,\qquad
N_2=n_2+n_3,\qquad
N_3=n_3,
\]
we have
\begin{align*}
N_1^2+N_2^2+N_3^2-3N_1-2N_2-N_3 &= n_1^2+2n_1n_2+2n_1n_3+2n_2^2+4n_2n_3+3n_3^2\\
&\qquad -3n_1-5n_2-6n_3.
\end{align*}

Thus, \eqref{eq:mainresult1} gives
\begin{align*}
&\sum_{n_1,n_2,n_3\geq0}
q^{n_1^2+2n_1n_2+2n_1n_3+2n_2^2+4n_2n_3+3n_3^2-3n_1-5n_2-6n_3}\\
&\quad\times
\begin{bmatrix}
n_1+10-2N_1\\
n_1
\end{bmatrix}_q
\begin{bmatrix}
n_2+17-2N_1-2N_2\\
n_2
\end{bmatrix}_q
\begin{bmatrix}
n_3+21-2N_1-2N_2-2N_3\\
n_3
\end{bmatrix}_q\\
&=
\sum_{s'=1}^{4}
G_{s'}((2,2,3),q)
G_{s'}((1,1,1),1/q).
\end{align*}

Here,
\begin{align*}
G_s((2,2,3),q)
={}&
\sum_{n_1,n_2,n_3\geq0}
q^{N_1^2+N_2^2+N_3^2+N_s+\cdots+N_3}\\
&\quad\times
\begin{bmatrix}
n_1+7-2N_1-\alpha_{1,s-1}\\
n_1
\end{bmatrix}_q
\begin{bmatrix}
n_2+12-2N_1-2N_2-\alpha_{2,s-1}\\
n_2
\end{bmatrix}_q\\
&\quad\times
\begin{bmatrix}
n_3+15-2N_1-2N_2-2N_3-\alpha_{3,s-1}\\
n_3
\end{bmatrix}_q.
\end{align*}
Similarly,
\begin{align*}
G_s((1,1,1),q)
={}&
\sum_{n_1,n_2,n_3\geq0}
q^{N_1^2+N_2^2+N_3^2+N_s+\cdots+N_3}\\
&\quad\times
\begin{bmatrix}
n_1+3-2N_1-\alpha_{1,s-1}\\
n_1
\end{bmatrix}_q
\begin{bmatrix}
n_2+5-2N_1-2N_2-\alpha_{2,s-1}\\
n_2
\end{bmatrix}_q\\
&\quad\times
\begin{bmatrix}
n_3+6-2N_1-2N_2-2N_3-\alpha_{3,s-1}\\
n_3
\end{bmatrix}_q.
\end{align*}



Substituting the relevant parameters into the definition of $G_s$ gives the following eight polynomials:
\begin{align*}
G_1((2,2,3),q)
={}&1+q^2+q^3+2q^4+2q^5+4q^6+3q^7+5q^8+5q^9\\
&\quad+7q^{10}+6q^{11}+9q^{12}+6q^{13}+8q^{14}+6q^{15}\\
&\quad+7q^{16}+3q^{17}+4q^{18}+q^{20},\\[4pt]
G_2((2,2,3),q)
={}&1+q+q^2+2q^3+3q^4+4q^5+6q^6+6q^7+8q^8\\
&\quad+10q^9+12q^{10}+12q^{11}+14q^{12}+14q^{13}+14q^{14}\\
&\quad+13q^{15}+12q^{16}+8q^{17}+7q^{18}+3q^{19}+q^{20},\\[4pt]
G_3((2,2,3),q)
={}&1+q+2q^2+2q^3+4q^4+5q^5+7q^6+8q^7+11q^8\\
&\quad+12q^9+16q^{10}+16q^{11}+19q^{12}+18q^{13}+20q^{14}\\
&\quad+17q^{15}+17q^{16}+12q^{17}+10q^{18}+4q^{19}+3q^{20},\\[4pt]
G_4((2,2,3),q)
={}&1+q+2q^2+3q^3+4q^4+5q^5+8q^6+9q^7+12q^8\\
&\quad+14q^9+17q^{10}+18q^{11}+22q^{12}+21q^{13}+22q^{14}\\
&\quad+20q^{15}+19q^{16}+14q^{17}+12q^{18}+6q^{19}+3q^{20}.
\end{align*}

Similarly,
\begin{align*}
G_1((1,1,1),q)&=1+q^2,\\
G_2((1,1,1),q)&=1+q+q^2+q^3,\\
G_3((1,1,1),q)&=1+q+2q^2+q^3,\\
G_4((1,1,1),q)&=1+q+2q^2+2q^3.
\end{align*}

Evaluating these polynomials explicitly and simplifying yields
\begin{align*}
&\sum_{n_1,n_2,n_3\geq0}
q^{n_1^2+2n_1n_2+2n_1n_3+2n_2^2+4n_2n_3+3n_3^2-3n_1-5n_2-6n_3}\\
&\quad\times
\begin{bmatrix}
n_1+10-2N_1\\
n_1
\end{bmatrix}_q
\begin{bmatrix}
n_2+17-2N_1-2N_2\\
n_2
\end{bmatrix}_q
\begin{bmatrix}
n_3+21-2N_1-2N_2-2N_3\\
n_3
\end{bmatrix}_q\\
&= 4q^{-3}+10q^{-2}+15q^{-1}+27+36q+53q^2+74q^3+99q^4+123q^5+157q^6+186q^7+221q^8\\&\qquad+249q^9+277q^{10}+283q^{11}+293q^{12}+275q^{13}+253q^{14}+208q^{15}+163q^{16}+99q^{17}+60q^{18}\\&\qquad+20q^{19}+8q^{20}.
\numberthis\label{eq:example23.3}
\end{align*}

Maple verifies \eqref{eq:example23.3}.

\section{New $\Vec{M}$ Versions of Andrews-Gordon Identities}\label{s4}
We now pass on to the new $\Vec{M}$ versions of the Andrews-Gordon Identities. Assume $M_1\ge M_2\ge M_3\ge\cdots\ge M_\nu(\ge 0)$. Then
\begin{align*}
&
\sum_{\mathbf n}
\frac{
q^{
N_1^2+\cdots+N_\nu^2
-
M_1N_1
-
M_2N_2
-
\cdots
-
M_\nu N_\nu
}
}{
(q)_{n_1}(q)_{n_2}\cdots(q)_{n_\nu}
}
\\
&=
\sum_{s'=1}^{\nu+1}
\prod_{\substack{
j\ge 1\\
j\not\equiv 0,\pm{s'}\!\!\!\!\pmod{2\nu+3}
}}
\frac1{(1-q^j)}
\,
I_{s'}
\!\left(
M_1-M_2,
M_2-M_3,
\dots,
M_\nu;
1/q
\right),\numberthis\label{eq:finiteMversion}
\end{align*}
where
\[
I_{s'}(M_1-M_2,M_2-M_3,\dots,M_\nu;q)
:=
\sum_{j=-\infty}^{\infty}
\Biggl(
q^{\,
j\bigl((2j+1)(2\nu+3)-2s'\bigr)}
\left[
\begin{matrix}
M_1-M_2,\,
M_2-M_3,\,
\dots,\,
M_\nu
\\[2mm]
(2\nu+3)j+\dfrac{\nu+1+a-s'}{2}
\end{matrix}
\right]_q
\]
\[
-
q^{\,
(2j+1)\bigl((2\nu+3)j+s'\bigr)}
\left[
\begin{matrix}
M_1-M_2,\,
M_2-M_3,\,
\dots,\,
M_\nu
\\[2mm]
(2\nu+3)j+\dfrac{\nu+1+a+s'}{2}
\end{matrix}
\right]_q
\Biggr),
\]
where
\[
a\equiv
\left(
\sum_{i=1}^{\nu}M_i+\nu+1+s'
\right)
\pmod{2},
\qquad
a\in\{0,1\}.
\]

It is easy to verify that \eqref{eq:finiteMversion} follows from \eqref{eq:mainresult1} in the limit $l_\nu\rightarrow\infty$. Now, we state two examples to illustrate \eqref{eq:finiteMversion}.

\subsection{Example 1}
Consider
$\nu=2$ with
\[
M_1=2,\qquad M_2=1.
\]
Then
\[
M_1-M_2=1,\qquad M_2=1,
\]
and $2\nu+3=7$. Moreover,
\[
a\equiv M_1+M_2+\nu+1+s'
   \equiv 2+1+2+1+s'
   \equiv s' \pmod 2.
\]
Thus $a=1$ for $s'=1,3$, while $a=0$ for $s'=2$.

Consequently, \eqref{eq:finiteMversion} gives
\[
\sum\limits_{n_1,n_2\geq0}
\frac{q^{N_1^2+N_2^2-2N_1-N_2}}
     {(q)_{n_1}(q)_{n_2}}
=
\sum\limits_{s'=1}^{3}
\prod\limits_{\substack{j\geq1\\j\not\equiv0,\pm s'\pmod 7}}
\frac{1}{1-q^j}\,
I_{s'}(1,1;1/q),
\]
where
\[
N_1=n_1+n_2,\qquad N_2=n_2,
\]
so that
\begin{align*}
N_1^2+N_2^2-2N_1-N_2 &= (n_1+n_2)^2+n_2^2-2(n_1+n_2)-n_2\\
&= n_1^2+2n_1n_2+2n_2^2-2n_1-3n_2.    
\end{align*}

The three polynomials occurring on the right-hand side are obtained by substituting $(M_1-M_2,M_2)=(1,1)$ into the definition of $I_{s'}$. Explicitly,
\[
\begin{aligned}
I_1(1,1;q)
&=\sum_{j=-\infty}^{\infty}
q^{j(14j+5)}
\left[\begin{matrix}1,1\\7j+\frac{3}{2}\end{matrix}\right]_q
-
\sum_{j=-\infty}^{\infty}
q^{(2j+1)(7j+1)}
\left[\begin{matrix}1,1\\7j+\frac{5}{2}\end{matrix}\right]_q\\
&=\sum_{j=-\infty}^{\infty}
q^{j(14j+5)}
\sum_{\substack{j_1,j_2\\j_1+j_2=7j+3}}
q^{j_1(1-j_2)}
\begin{bmatrix}1\\j_2\end{bmatrix}_q
\begin{bmatrix}1+j_2\\j_1\end{bmatrix}_q\\
&\qquad-\sum_{j=-\infty}^{\infty}
q^{(2j+1)(7j+1)}
\sum_{\substack{j_1,j_2\\j_1+j_2=7j+4}}
q^{j_1(1-j_2)}
\begin{bmatrix}1\\j_2\end{bmatrix}_q
\begin{bmatrix}1+j_2\\j_1\end{bmatrix}_q\\
&= 1,
\end{aligned}
\]
\[
\begin{aligned}
I_2(1,1;q)
&=\sum_{j=-\infty}^{\infty}
q^{j(14j+3)}
\left[\begin{matrix}1,1\\7j+\frac{1}{2}\end{matrix}\right]_q
-
\sum_{j=-\infty}^{\infty}
q^{(2j+1)(7j+2)}
\left[\begin{matrix}1,1\\7j+\frac{5}{2}\end{matrix}\right]_q\\
&=\sum_{j=-\infty}^{\infty}
q^{j(14j+3)}
\sum_{\substack{j_1,j_2\\j_1+j_2=7j+2}}
q^{j_1(1-j_2)}
\begin{bmatrix}1\\j_2\end{bmatrix}_q
\begin{bmatrix}1+j_2\\j_1\end{bmatrix}_q\\
&\qquad-\sum_{j=-\infty}^{\infty}
q^{(2j+1)(7j+2)}
\sum_{\substack{j_1,j_2\\j_1+j_2=7j+4}}
q^{j_1(1-j_2)}
\begin{bmatrix}1\\j_2\end{bmatrix}_q
\begin{bmatrix}1+j_2\\j_1\end{bmatrix}_q\\
&= 1+q,
\end{aligned}
\]
and
\[
\begin{aligned}
I_3(1,1;q)
&=\sum_{j=-\infty}^{\infty}
q^{j(14j+1)}
\left[\begin{matrix}1,1\\7j+\frac{1}{2}\end{matrix}\right]_q
-
\sum_{j=-\infty}^{\infty}
q^{(2j+1)(7j+3)}
\left[\begin{matrix}1,1\\7j+\frac{7}{2}\end{matrix}\right]_q\\
&=\sum_{j=-\infty}^{\infty}
q^{j(14j+1)}
\sum_{\substack{j_1,j_2\\j_1+j_2=7j+2}}
q^{j_1(1-j_2)}
\begin{bmatrix}1\\j_2\end{bmatrix}_q
\begin{bmatrix}1+j_2\\j_1\end{bmatrix}_q\\
&\qquad-\sum_{j=-\infty}^{\infty}
q^{(2j+1)(7j+3)}
\sum_{\substack{j_1,j_2\\j_1+j_2=7j+5}}
q^{j_1(1-j_2)}
\begin{bmatrix}1\\j_2\end{bmatrix}_q
\begin{bmatrix}1+j_2\\j_1\end{bmatrix}_q\\
&= 1+q.
\end{aligned}
\]

Therefore, the identity in this particular case takes the form
\begin{align*}
\sum\limits_{n_1,n_2\geq0}\frac{q^{n_1^2+2n_1n_2+2n_2^2-2n_1-3n_2}}{(q)_{n_1}(q)_{n_2}} &= \prod\limits_{\substack{j\geq1\\j\not\equiv0,\pm1\pmod7}}\frac{1}{1-q^j} + \prod\limits_{\substack{j\geq1\\j\not\equiv0,\pm2\pmod7}}\frac{1+q^{-1}}{1-q^j}\\
&\qquad + \prod\limits_{\substack{j\geq1\\j\not\equiv0,\pm3\pmod7}}\frac{1+q^{-1}}{1-q^j}.\numberthis\label{eq:example1}
\end{align*}
This provides a concrete instance with $\nu = 2$ of the $\vec M$ version of the Andrews-Gordon identity. Maple verifies \eqref{eq:example1}.

\subsection{Example 2}
We now consider
$\nu=2$ with
\[
M_1=3,\qquad M_2=2.
\]
Then
\[
M_1-M_2=1,\qquad M_2=2,
\]
and $2\nu+3=7$. Moreover,
\[
a\equiv M_1+M_2+\nu+1+s'
   \equiv 3+2+2+1+s'
   \equiv s' \pmod 2.
\]
Thus $a=1$ for $s'=1,3$, while $a=0$ for $s'=2$.

Consequently, \eqref{eq:finiteMversion} gives
\[
\sum\limits_{n_1,n_2\geq0}
\frac{q^{N_1^2+N_2^2-3N_1-2N_2}}
     {(q)_{n_1}(q)_{n_2}}
=
\sum\limits_{s'=1}^{3}
\prod\limits_{\substack{j\geq1\\j\not\equiv0,\pm s'\pmod 7}}
\frac{1}{1-q^j}\,
I_{s'}(1,1;1/q),
\]
where
\[
N_1=n_1+n_2,\qquad N_2=n_2,
\]
so that
\begin{align*}
N_1^2+N_2^2-3N_1-2N_2 &= (n_1+n_2)^2+n_2^2-3(n_1+n_2)-2n_2\\
&= n_1^2+2n_1n_2+2n_2^2-3n_1-5n_2.    
\end{align*}

The three polynomials occurring on the right-hand side are obtained by substituting $(M_1-M_2,M_2)=(1,2)$ into the definition of $I_{s'}$. Explicitly,
\[
\begin{aligned}
I_1(1,2;q)
&=\sum_{j=-\infty}^{\infty}
q^{j(14j+5)}
\left[\begin{matrix}1,2\\7j+\frac{3}{2}\end{matrix}\right]_q
-
\sum_{j=-\infty}^{\infty}
q^{(2j+1)(7j+1)}
\left[\begin{matrix}1,2\\7j+\frac{5}{2}\end{matrix}\right]_q\\
&=\sum_{j=-\infty}^{\infty}
q^{j(14j+5)}
\sum_{\substack{j_1,j_2\\j_1+j_2=7j+4}}
q^{j_1(2-j_2)}
\begin{bmatrix}2\\j_2\end{bmatrix}_q
\begin{bmatrix}1+j_2\\j_1\end{bmatrix}_q\\
&\qquad-\sum_{j=-\infty}^{\infty}
q^{(2j+1)(7j+1)}
\sum_{\substack{j_1,j_2\\j_1+j_2=7j+5}}
q^{j_1(2-j_2)}
\begin{bmatrix}2\\j_2\end{bmatrix}_q
\begin{bmatrix}1+j_2\\j_1\end{bmatrix}_q\\
&= 1+q^2,
\end{aligned}
\]
\[
\begin{aligned}
I_2(1,2;q)
&=\sum_{j=-\infty}^{\infty}
q^{j(14j+3)}
\left[\begin{matrix}1,2\\7j+\frac{1}{2}\end{matrix}\right]_q
-
\sum_{j=-\infty}^{\infty}
q^{(2j+1)(7j+2)}
\left[\begin{matrix}1,2\\7j+\frac{5}{2}\end{matrix}\right]_q\\
&=\sum_{j=-\infty}^{\infty}
q^{j(14j+3)}
\sum_{\substack{j_1,j_2\\j_1+j_2=7j+3}}
q^{j_1(2-j_2)}
\begin{bmatrix}2\\j_2\end{bmatrix}_q
\begin{bmatrix}1+j_2\\j_1\end{bmatrix}_q\\
&\qquad-\sum_{j=-\infty}^{\infty}
q^{(2j+1)(7j+2)}
\sum_{\substack{j_1,j_2\\j_1+j_2=7j+5}}
q^{j_1(2-j_2)}
\begin{bmatrix}2\\j_2\end{bmatrix}_q
\begin{bmatrix}1+j_2\\j_1\end{bmatrix}_q\\
&= 1+q+q^2+q^3,
\end{aligned}
\]
and
\[
\begin{aligned}
I_3(1,2;q)
&=\sum_{j=-\infty}^{\infty}
q^{j(14j+1)}
\left[\begin{matrix}1,2\\7j+\frac{1}{2}\end{matrix}\right]_q
-
\sum_{j=-\infty}^{\infty}
q^{(2j+1)(7j+3)}
\left[\begin{matrix}1,2\\7j+\frac{7}{2}\end{matrix}\right]_q\\
&=\sum_{j=-\infty}^{\infty}
q^{j(14j+1)}
\sum_{\substack{j_1,j_2\\j_1+j_2=7j+3}}
q^{j_1(2-j_2)}
\begin{bmatrix}2\\j_2\end{bmatrix}_q
\begin{bmatrix}1+j_2\\j_1\end{bmatrix}_q\\
&\qquad-\sum_{j=-\infty}^{\infty}
q^{(2j+1)(7j+3)}
\sum_{\substack{j_1,j_2\\j_1+j_2=7j+6}}
q^{j_1(2-j_2)}
\begin{bmatrix}2\\j_2\end{bmatrix}_q
\begin{bmatrix}1+j_2\\j_1\end{bmatrix}_q\\
&= 1+q+2q^2+q^3.
\end{aligned}
\]

Therefore, the identity in this particular case takes the form
\begin{align*}
\sum\limits_{n_1,n_2\geq0}\frac{q^{n_1^2+2n_1n_2+2n_2^2-3n_1-5n_2}}{(q)_{n_1}(q)_{n_2}} &= \prod\limits_{\substack{j\geq1\\j\not\equiv0,\pm1\pmod7}}\frac{1+q^{-2}}{1-q^j} + \prod\limits_{\substack{j\geq1\\j\not\equiv0,\pm2\pmod7}}\frac{1+q^{-1}+q^{-2}+q^{-3}}{1-q^j}\\
&\qquad + \prod\limits_{\substack{j\geq1\\j\not\equiv0,\pm3\pmod7}}\frac{1+q^{-1}+2q^{-2}+q^{-3}}{1-q^j}.\numberthis\label{eq:example2}
\end{align*}
This provides another instance with $\nu = 2$ of the $\vec{M}$ version of the Andrews-Gordon identity. Maple verifies \eqref{eq:example2}.

\subsection{Example 3}
We now consider $\nu=3$ with
\[
M_1=3,\qquad M_2=2,\qquad M_3=1.
\]
Then
\[
M_1-M_2=M_2-M_3=M_3=1,
\]
and
\[
2\nu+3=9.
\]
Moreover,
\[
a\equiv M_1+M_2+M_3+\nu+1+s'
\equiv 3+2+1+3+1+s'
\equiv s' \pmod 2.
\]
Thus $a=1$ for $s'=1,3$, while $a=0$ for $s'=2,4$.

Consequently, \eqref{eq:finiteMversion} gives
\[
\sum_{n_1,n_2,n_3\geq0}
\frac{q^{N_1^2+N_2^2+N_3^2-3N_1-2N_2-N_3}}
     {(q)_{n_1}(q)_{n_2}(q)_{n_3}}
=
\sum_{s'=1}^{4}
\prod_{\substack{j\geq1\\j\not\equiv0,\pm s'\pmod9}}
\frac{1}{1-q^j}
I_{s'}(1,1,1;1/q),
\]
where
\[
N_1=n_1+n_2+n_3,\qquad
N_2=n_2+n_3,\qquad
N_3=n_3.
\]
Hence
\[
\begin{aligned}
N_1^2+N_2^2+N_3^2-3N_1-2N_2-N_3
={}&n_1^2+2n_1n_2+2n_1n_3-3n_1\\
&+2n_2^2+4n_2n_3-5n_2
+3n_3^2-6n_3.
\end{aligned}
\]

The four polynomials occurring on the right-hand side are obtained by substituting $(M_1-M_2,M_2-M_3,M_3)=(1,1,1)$ into the definition of $I_{s'}$. Explicitly,
\[
\begin{aligned}
I_1(1,1,1;q)&=1+q^2,\\
I_2(1,1,1;q)&=1+q+q^2+q^3,\\
I_3(1,1,1;q)&=1+q+2q^2+q^3,\\
I_4(1,1,1;q)&=1+q+2q^2+2q^3.
\end{aligned}
\]

Therefore,
\[
\begin{aligned}
I_1(1,1,1;1/q)&=1+q^{-2},\\
I_2(1,1,1;1/q)&=1+q^{-1}+q^{-2}+q^{-3},\\
I_3(1,1,1;1/q)&=1+q^{-1}+2q^{-2}+q^{-3},\\
I_4(1,1,1;1/q)&=1+q^{-1}+2q^{-2}+2q^{-3}.
\end{aligned}
\]

Thus, the identity in this particular case takes the form
\begin{align*}
&\sum_{n_1,n_2,n_3\geq0}\frac{q^{n_1^2+2n_1n_2+2n_1n_3-3n_1+2n_2^2+4n_2n_3-5n_2+3n_3^2-6n_3}}{(q)_{n_1}(q)_{n_2}(q)_{n_3}}\\
&\qquad = \prod\limits_{\substack{j\geq1\\j\not\equiv0,\pm1\pmod9}}\frac{1+q^{-2}}{1-q^j}+\prod\limits_{\substack{j\geq1\\j\not\equiv0,\pm2\pmod9}}\frac{1+q^{-1}+q^{-2}+q^{-3}}{1-q^j}\\
&\qquad\quad + \prod\limits_{\substack{j\geq1\\j\not\equiv0,\pm3\pmod9}}\frac{1+q^{-1}+2q^{-2}+q^{-3}}{1-q^j}+\prod\limits_{\substack{j\geq1\\j\not\equiv0,\pm4\pmod9}}\frac{1+q^{-1}+2q^{-2}+2q^{-3}}{1-q^j}.\numberthis\label{eq:example3}
\end{align*}

This provides an instance with $\nu=3$ of the $\vec M$ version of the Andrews-Gordon identity. Maple verifies \eqref{eq:example3}.

\section*{Acknowledgements}
We thank Ae Ja Yee and Jihyeug Jang for their curiosity about \cite[eq. (5.7)]{BerkovichPaule2} which reignited our interest in this work.

\end{document}